\documentclass[letterpaper]{amsart}

\usepackage[letterpaper,left=1.5in,right=1.5in,top=1in,bottom=1in]{geometry}

\title[Serre-Godeaux varieties and the \'etale index]{Serre-Godeaux varieties and the \'etale index}

\author[Benjamin~Antieau]{Benjamin Antieau$^{*}$}
\address{UCLA, Department of Mathematics, 520 Portola Plaza, Los Angeles, CA 90095-1555}
\email{antieau@math.ucla.edu}
\thanks{$^{*}$The author was supported in part by the NSF under Grant RTG DMS 0838697}

\author[Ben~Williams]{Ben Williams}
\address{USC, Department of Mathematics, 3620 South Vermont Avenue, Los Angeles, CA 90089-2532}
\email{tbwillia@usc.edu}


\newcommand{\myauthor}{Benjamin Antieau and Ben Williams}
\newcommand{\mytitle}{Serre-Godeaux varieties and the \'etale index}

\usepackage[pdfstartview=FitH,
            pdfauthor={\myauthor},
            pdftitle={\mytitle},
            colorlinks,
            linkcolor=reference,
            citecolor=citation,
            urlcolor=e-mail,
            ]{hyperref}

\usepackage{amsmath}
\usepackage{amscd}
\usepackage{amsbsy}
\usepackage{amssymb}
\usepackage{verbatim}
\usepackage{eufrak}
\usepackage{microtype}
\usepackage{hyperref}
\usepackage{mathrsfs}
\usepackage{amsthm}
\usepackage[all,cmtip]{xy}
\usepackage{tikz}
\usetikzlibrary{matrix,arrows}
\usepackage[abbrev,lite]{amsrefs}

\usepackage{color}
\definecolor{todo}{rgb}{1,0,0}
\definecolor{conditional}{rgb}{0,1,0}
\definecolor{e-mail}{rgb}{0,.40,.80}
\definecolor{reference}{rgb}{.20,.60,.22}
\definecolor{mrnumber}{rgb}{.80,.40,0}
\definecolor{citation}{rgb}{0,.40,.80}

\DeclareMathAlphabet{\mathpzc}{OT1}{pzc}{m}{it}
\usepackage{dsfont}

\DeclareMathOperator{\PGL}{PGL}

\DeclareMathOperator{\GL}{GL}

\DeclareMathOperator{\tors}{tors}

\DeclareMathOperator{\un}{un}

\DeclareMathOperator{\Spec}{Spec}

\DeclareMathOperator{\Hoh}{H}
\DeclareMathOperator{\Eoh}{E}

\DeclareMathOperator{\Br}{Br}

\DeclareMathOperator*{\holim}{holim}

\DeclareMathOperator{\eti}{eti}
\DeclareMathOperator{\ind}{ind}

\DeclareMathOperator{\per}{per}

\newcommand{\topo}{{\mathrm{top}}}

\newcommand{\et}{\mathrm{\acute{e}t}}

\DeclareFontEncoding{OT2}{}{}

\newcommand{\K}{\mathbf{K}}

\DeclareMathOperator{\KU}{KU}

\newcommand{\iso}{\cong}

\newcommand{\CC}{\mathds{C}}

\newcommand{\QQ}{\mathds{Q}}
\newcommand{\ZZ}{\mathds{Z}}

\newcommand{\PP}{\mathds{P}}

\newcommand{\Gm}{\mathds{G}_{m}}

\let\oldmarginpar\marginpar
\renewcommand\marginpar[1]{\-\oldmarginpar[\raggedleft\footnotesize #1]%
{\raggedright\footnotesize #1}}

\theoremstyle{plain}
\newtheorem{theorem}{Theorem}[section]

\newtheorem{scholium}[theorem]{Scholium}
\newtheorem{proposition}[theorem]{Proposition}

\newtheorem{corollary}[theorem]{Corollary}

\theoremstyle{definition}

\theoremstyle{remark}
\newtheorem{remark}[theorem]{Remark}

\begin{document}

\begin{abstract}
    We use the Serre-Godeaux varieties of finite groups, projective representation theory, the twisted Atiyah-Segal
    completion theorem, and our
    previous work on the topological period-index problem to
    compute the etale index of Brauer classes $\alpha\in\Br(X)$ in some specific examples.
    In particular, these computations show that
    the \'etale index of $\alpha$ differs from the period of $\alpha$ in general. As an application, we
    compute the index of unramified classes in the function fields of high-dimensional
    Serre-Godeaux varieties in terms of projective representation theory.
\end{abstract}

\maketitle

Let $X$ be a connected scheme, and let $\Br'(X)=\Hoh^2_{\et}(X,\Gm)_{\tors}$ be its cohomological Brauer
group. There is a subgroup $\Br(X)\subseteq\Br'(X)$ consisting of those Brauer classes that are
represented by an Azumaya algebra. A class $\alpha\in\Br'(X)$ is in $\Br(X)$ if and only if it is in
the image of the coboundary map
\[
    \Hoh^1_{\et}(X,\PGL_n)\xrightarrow{\delta_n}\Hoh^2_{\et}(X,\Gm)
\]
for some $n$, where the coboundary arises from the central extension
\[
	1\rightarrow\Gm\rightarrow\GL_n\rightarrow\PGL_n\rightarrow 1.
\]
The pointed cohomology set $\Hoh^1_{\et}(X,\PGL_n)$ classifies $\PGL_n$-torsors, $\PP^{n-1}$-bundles, or degree-$n$ Azumaya algebras.
Recall that a degree-$n$ Azumaya algebra is an algebra on $X$ which looks \'etale-locally like the algebra of $n\times n$ matrices over the ring of regular functions on $X$.
The class $\delta_n(P)$ of a $\PGL_n$-torsor $P$ is precisely the obstruction to lifting $P$ to a
$\GL_n$-torsor. Viewed from the perspective of Azumaya algebras, this class is the obstruction to writing an Azumaya
algebra $A$ as the endomorphism algebra of a vector bundle. The index of $\alpha\in\Br(X)$ is
\[
    \ind(\alpha)=\gcd\{n|\text{$\alpha$ is in the image of $\delta_n$\}}.
\]
When $X=\Spec k$ is the spectrum of a field, the index of a class $\alpha\in\Br(k)$ is the degree of the unique division algebra representing $\alpha$.

The period of $\alpha$, denoted $\per(\alpha)$, is the order of $\alpha\in\Br'(X)$. In general,
\[
    \per(\alpha)|\ind(\alpha).
\]
When $X=\Spec k$, the two integers have the same prime divisors.
It is a major open problem in the theory of Azumaya algebras, even when $X=\Spec k$, to determine the possible pairs $(\per(\alpha),\ind(\alpha))$ for $\alpha\in\Br(X)$ where $X$ is
fixed, or where the dimension of $X$ is fixed. For an introduction to what is known and for
further references, see~\cite{aw1}*{Section~1}.

In~\cite{antieau}, the first author constructed an invariant $\eti(\alpha)$ of cohomological Brauer classes
$\alpha$, called the \'etale index. It is constructed as the positive generator of a rank map
$\K^{\alpha,\et}_0(X)\rightarrow\ZZ$ from twisted \'etale $K$-theory, which was first introduced in~\cite{antieau}.
It was proved there and in~\cite{antieau-cech} that, for $\alpha\in\Br(X)$,
\begin{equation}\label{eq:div1}
    \per(\alpha)|\eti(\alpha)|\ind(\alpha),
\end{equation}
and that, in general, the \'etale index is strictly smaller than the index. The question of whether
or not the \'etale index ever differed from the period was left open. In this paper, we
present in Theorem~\ref{thm:main} a class of examples arising from finite group theory for which the \'etale index is
indeed different from the period. The examples are
Serre-Godeaux varieties~\cite{serre}, which are smooth projective varieties that approximate
the homotopy type of $BG\times K(\ZZ,2)$. Our approach to these
varieties comes from the presentation in the proof of~\cite{atiyah-hirzebruch}*{Proposition~6.6}.

The \'etale index is interesting for at least the following two reasons: first it has the property that it is finite for all classes
$\alpha\in\Br'(X)$, provided $X$ itself is of finite \'etale cohomological dimension, this contrasts
with the index \textit{per se} which may not be defined when $\Br(X)\neq\Br'(X)$. Second, in~\cite{antieau} we gave upper bounds
for the \'etale index in terms of the period using stable homotopy theory, thus solving an analogue of the still open period-index
conjecture.

If $X$ is a smooth projective complex variety, the unramified Brauer group of the function
field $\CC(X)$ is $\Br_{\un}(\CC(X))=\Br(X)$. If $G$ is a finite group and if $X$ is a Serre-Godeaux variety
of sufficiently high dimension, our method allows us to compute in
Theorem~\ref{thm:indsg} the
index of $\alpha$ for $\alpha\in\Br_{\un}(X)$ in terms of projective representation theory.
This computation does not appear accessible by other methods.

\subsection*{Acknowledgments}
We would like to thank Anssi Lahtinen for several very useful exchanges about the
completion theorem in twisted $K$-theory~\cite{lahtinen}.

\section{Definitions}

When $X$ is a topological space, the cohomological Brauer group is
$\Br'_{\topo}(X)=\Hoh^3(X,\ZZ)_{\tors}$, and the Brauer group $\Br_\topo(X)\subseteq\Br_\topo'(X)$ is the
subgroup of all classes represented by a topological Azumaya algebra.
If $X$ is a finite CW-complex, then $\Br_{\topo}(X)=\Br'_{\topo}(X)$ by an argument of Serre~\cite{grothendieck-brauer}.
For a class $\alpha\in\Br'_{\topo}(X)$, we define $\per_{\topo}(\alpha)$ to be the order of $\alpha$.
There are coboundary maps
\begin{equation*}
    \Hoh^1(X,\PGL_n(\CC))\xrightarrow{\delta_n}\Hoh^2(X,\CC^*)\iso\Hoh^3(X,\ZZ)_{\tors},
\end{equation*}
and we define
\[
    \ind_{\topo}(\alpha)=\gcd\{n|\text{$\alpha$ is in the image of $\delta_n$\}}
\]
when $\alpha\in\Br_{\topo}(X)$.
There is an $\alpha$-twisted $K$-theory spectrum, $\KU(X)_{\alpha}$, as defined in~\cite{atiyah-segal}. If $X$ is a
finite CW-complex the group $\KU^0(X)_{\alpha}$ is the Grothendieck of $\alpha$-twisted vector bundles.
This group is equipped with a natural rank map $\KU^0(X)_{\alpha}\rightarrow\ZZ$. The $K$-theoretic index
$\ind_{\K}(\alpha)$ is defined in~\cite{aw1} to be the positive generator of the image of the rank map. In topology, the $K$-theoretic index plays
the same role as the \'etale index $\eti(\alpha)$ does in algebraic geometry.
In general, we have the following analog of~\eqref{eq:div1}:
\[
\per_{\topo}(\alpha)|\ind_{\K}(\alpha)|\ind_{\topo}(\alpha).
\]
In the case where $X$ is a finite CW-complex, because topological twisted
$K$-theory satisfies descent, we showed in~\cite{aw1}*{Lemma~2.23} that $\ind_{\K}(\alpha)=\ind_{\topo}(\alpha)$. 

If $X$ is a complex algebraic scheme, there is a natural map
$\Br'(X)\rightarrow\Hoh^3(X,\ZZ)_{\tors}$.
If $\alpha\in\Br'(X)$ is a cohomological Brauer class, we write $\overline{\alpha}$ for the image
of $\alpha$ in $\Hoh^3(X,\ZZ)_{\tors}$. In general, $\per_{\topo}(\overline{\alpha})$ divides $\per(\alpha)$, and
it is easy to see that $\ind_{\topo}(\overline{\alpha})$ divides $\ind(\alpha)$ as well, because any algebraic Azumaya algebra gives rise to a topological Azumaya algebra. More importantly
for the purposes of this paper, a result
of~\cite{aw1} says that if $X$ has the homotopy type of a finite CW-complex, then
\begin{equation} \label{eq:1}
    \ind_{\topo}(\overline{\alpha}) =\ind_\K(\overline \alpha) |\eti(\alpha).
\end{equation}
Fundamentally, this result again depends on descent for topological twisted $K$-theory.

Let $X$ be a topological space and $\alpha\in\Br_{\topo}'(X)$.
There is a twisted Atiyah-Hirzebruch spectral sequence~\cite{atiyah-segal-cohomology}
\[
    \Eoh_2^{p,q}=\Hoh^p\left(X,\KU^q(*)\right)\Longrightarrow\KU^{p+q}(X)_{\alpha},
\]
where $\KU^*(X)_{\alpha}$ denotes the $\alpha$-twisted $K$-theory of $X$. The differentials $d_r^{\alpha}$ are of degree
$(r,1-r)$, and $d_r^{\alpha}=0$ unless $r$ is odd.
The first non-trivial differential $d_3^{\alpha}:\Hoh^0(X,\ZZ)\rightarrow\Hoh^3(X,\ZZ)$ has the property that $d_3^{\alpha}(1)=\pm\alpha\in\Hoh^3(X,\ZZ)$, by
\cite{atiyah-segal} or~\cite{antieau}. This spectral sequence is natural, in the
sense that if $f:X\rightarrow Y$ is a map and $\alpha\in\Hoh^3(Y,\ZZ)_{\tors}$, then there is a morphism from
the spectral sequence computing $\KU^*(Y)_{\alpha}$ to the spectral sequence computing
$\KU^*(X)_{f^*\alpha}$.

The twisted Atiyah-Segal spectral sequence converges strongly when $X$ is a finite
CW-complex. When $X$ is connected, the rank map $\KU^0(X)_{\alpha}\rightarrow\ZZ$ may be identified with the edge map
$\KU^0(X)_{\alpha}\rightarrow\Hoh^0(X,\ZZ)\iso\ZZ$ in the spectral
sequence, and the image of the rank map is then the group of
permanent cycles $\Eoh_\infty^{0,0}\subseteq\Hoh^0(X,\ZZ)$. The kernel of the rank map we
denote by $\widetilde{\KU}^0(X)_{\alpha}$. In general the spectral sequence converges only conditionally,
but something may still be said about the permanent cycles in special cases.

\begin{proposition}\label{prop:kah}
    Let $X=\cup_k X_k$ be a skeletal filtration of the connected CW-complex $X$ by finite
    connected CW-complexes with $X_0=*$,
    and let $\alpha \in \Hoh^3(X, \ZZ)_{\tors}$. Suppose that
    \begin{equation*}
        {\lim}^1\widetilde{\KU}^0(X_k)_{f_k^*\alpha}=0,
    \end{equation*}
    where $f_k$ is the inclusion of $X_k$ in $X$. Then, $\ind_\K(\alpha)$ is
    the positive generator of the subgroup of permanent cycles in $\Eoh_2^{0,0} \iso \ZZ$.
    \begin{proof}
        The natural map $\KU(X)_{\alpha}\rightarrow\holim_k\KU(X_k)_{f_k^*\alpha}$ is a homotopy equivalence, which
        induces a commutative diagram (the commutativity of the top square is deduced from the naturality of the relation between $\Eoh_\infty$
        and the filtration on the target group in a conditionally convergent spectral sequence,
        for which see \cite{boardman}*{Chapters 5 \& 7})
        \begin{equation*}
        \xymatrix{\KU^0(X)_{\alpha} \ar[d] \ar@{->>}[r] &   \lim_k\KU^0(X_k)_{f_k^*\alpha} \ar^h@{->>}[d]  \ar[r] & \KU^0(X_0) \iso \ZZ\\
          \Eoh_{\infty}^{0,0}(X) \ar@{^(->}[d] \ar[r] &  \lim_k\Eoh_{\infty}^{0,0}(X_k) \ar[ur] \ar@{->}[d] \\
          \Eoh_2^{0,0}(X) \ar@{=}[r] & \lim_k \Eoh_2^{0,0}(X_k) \ar@{=}[uur]},
        \end{equation*}
        where the top-left horizontal arrow is surjective by an argument analogous to that
        of Milnor in the untwisted case~\cite{milnor}.

        The map $h$ is surjective because of the hypothesis that
         \begin{equation*}
            {\lim}^1\widetilde{\KU}^0(X_k)_{f_k^*\alpha}=0.
        \end{equation*}
        It follows by a diagram-chase that $\Eoh_\infty^{0,0}(X) \subset \Eoh_2^{0,0}(X)$ is precisely the image of the rank map.
        %   The groups $\Eoh_\infty^{0,0}(X_k)$ form a descending chain of subgroups of
        %   $\Eoh_2^{0,0}(X_0) = \ZZ$. If this chain of subgroups does not stabilize, then $\lim_k \Eoh_\infty^{0,0}(X_k) = 0$, and the result is
        %   trivially true. Suppose then that the chain of subgroups does stabilize, that is there exists some $k$ such that $\Eoh_\infty^{0,0}(X_r)
        %   \to \Eoh_\infty^{0,0}(X_k) \iso \ZZ$ is an isomorphism for $r \ge k$. Since $\ZZ$ is projective,
        %   we can find a splitting of $\KU^*(X_k)_{f_k^*\alpha} \to \Eoh_\infty^{0,0}(X_k) \iso \ZZ$. Suppose
        %   now that in the diagram
        %   \[ \xymatrix{ \KU^*(X_{r+1})_{f_{r+1}^* \alpha} \ar^{h_{r+1}}[d] \ar[r] & \KU^*(X_r)_{f_r^* \alpha} = A_r
        %     \oplus \Eoh_\infty^{0,0}(X_r) \ar[d] 
        %     \\ \Eoh_\infty^{0,0}(X_{r+1}) \ar@{=}[r] & \Eoh_\infty^{0,0}(X_r) \iso \ZZ \ar@/^1em/[u] }\]
        %   a splitting, as denoted by the dashed arrow, has already been constructed. We
        %     inductively find compatible
        %   splittings
        %   \begin{equation*}
        %     \xymatrix{ \dots \ar[r] & \KU^*(X_r)_{f^*_r(\alpha)} \ar@{->>}[d] \ar[r] & \dots \ar[r] &
        %       \KU^*(X_k)_{f^*_k(\alpha)} \ar@{->>}[d] \\ \dots \ar@{=}[r] & \ZZ  \ar@/^1em/@{..>}[u] \ar@{=}[r] & \dots \ar@{=}[r] &  \ZZ \ar@/^1em/@{..>}[u]. }
        %   \end{equation*}
        %   We now have a splitting $\lim_k \Eoh^{0,0}_\infty \to \lim_k \KU^*(X_k)_{f^*_k(\alpha)}$, and in particular a surjection.
        % 
    \end{proof}
\end{proposition}

When $G$ is a topological group, there is yet another index associated to a class $\alpha\in\Hoh^3(BG,\ZZ)_{\tors}$.
Given a projective representation $\pi:G\rightarrow\PGL_n(\CC)$, there is an associated class
$[\pi]$ in $\Hoh^1(BG,\PGL_n(\CC))$. 
We define the representation index
\[
\ind_G(\alpha)=\gcd\{n|\text{there exists $\pi:G\rightarrow\PGL_n$ with $\delta_n([\pi])=\alpha$}\}.
\]
% Note that not every map $BG\rightarrow B\U_n$ arises from a map $G\rightarrow \U_n$, even when $G$
% is finite. % TBJW: I'm commenting this out for now.
There are relations
\[
\per_{\topo}(\alpha)|\ind_{\K}(\alpha)|\ind_{\topo}(\alpha)|\ind_G(\alpha)
\]
for $\alpha\in\Hoh^3(BG,\ZZ)_{\tors}$. In the next section, we prove that these three indices
coincide when $G$ is a compact Lie group.

\section{Twisted equivariant K-theory}

% Let $G$ be a topological group, and let $X$ be a $G$-CW-complex. Set $X_G=(X\times EG)/G$.
% The equivariant cohomology group $\Hoh^3_G(X,\ZZ)_{\tors}=\Hoh^3(X_G,\ZZ)$ is the
% $G$-equivariant cohomological Brauer group $\Br_G'(X)$ of $X$. Just as in the
% non-equivariant theory, there is a subgroup, $\Br_G(X)$ consisting of those equivariant
% Brauer classes which are represnted by a $G$-equivariant Azumaya algebra, or equivalently,
% by a $G$-equivariant principle $\PGL_n$-bundle for some $n$.

If $G$ is a compact Lie group and $X$ is a finite $G$-CW-complex, we let
$\Hoh^n_G(X,\ZZ)=\Hoh^n(X_G,\ZZ)$, where $X_G=(X\times EG)/G$, and where $EG$ is a universal
$G$-space. When $\alpha\in\Hoh^3_G(X,\ZZ)$, there is a twisted equivariant $K$-theory
$\KU_G^*(X)_{\alpha}$ on $X$~\cite{atiyah-segal} and a twisted $K$-theory
$\KU^*(X_G)_{\alpha}$ on $X_G$. These are related by functorial morphisms
$\KU^*_G(X)_{\alpha}\rightarrow\KU^*(X_G)_{\alpha}$. The twisted equivariant $K$-theory is a module over
$\KU_G^0(\ast)\iso R(G)$. Let $I$ be the augmentation ideal of $R(G)$.
Below, $EG^n$ denotes the filtration on $EG$ given by Milnor~\cite{milnor-construction}, so that $EG^n$ is the join of
$(n+1)$ copies of $G$. The group $G$ acts freely on $EG^n$ with quotient space $BG^n$. The
spaces $BG^n$ form an exhaustive filtration of $BG$ by finite CW-complexes.

\begin{theorem}[Lahtinen~\cite{lahtinen}]
    If $G$ is a compact Lie group and
    $\alpha\in\Hoh^3(BG,\ZZ)=\Hoh^3_G(*,\ZZ)$, then there is a natural isomorphism of pro-groups
    \begin{equation*}
        \{\KU_G^0(*)_{\alpha}/I^n\cdot\KU_G^0(*)_{\alpha}\}\rightarrow\{\KU_G^0(EG^n)_{\alpha}\}.
    \end{equation*}
\end{theorem}

% \begin{remark}
%     Note that, in the twisted case, this morphism is really a morphism of pro-groups. There
%     might not be maps
%     \begin{equation*}
%         \KU_G^0(X)_{\alpha}/I^n\cdot\KU_G^0(X)_{\alpha}\rightarrow\KU_G^0(X\times EG^n)_{\alpha}.
%     \end{equation*}
%     However, for every $n$ there will exist some $m$ and a map
%     \begin{equation*}
%         \KU_G^0(X)_{\alpha}/I^n\cdot\KU_G^0(X)_{\alpha}\rightarrow\KU_G^0(X\times EG^m)_{\alpha}.
%     \end{equation*}
%     It is these which glue into the pro-group morphism. See~\cite{lahtinen}*{Lemma 5}.
% \end{remark}

\begin{corollary}\label{cor:asiso}
    In the situation above, the natural map
    \begin{equation*}
        \KU_G^0(*)_{\alpha}\hat{_I}\rightarrow\KU^0(BG)_{\alpha}
    \end{equation*}
    is an isomorphism, where $\KU_G^0(*)_{\alpha}\hat{_I}$ is the completion of
    $\KU^0_G(*)_{\alpha}$ along $I$.
\end{corollary}

Recall that a sequence $\cdots\rightarrow A_n\rightarrow A_{n-1}\rightarrow\cdots$ satisfies
the Mittag-Leffler condition if for every $k$, the image of $A_r\rightarrow A_k$ stabilizes
for $r$ sufficiently large. To be precise, for every $r$, there exists $n$ such that if
$k\geq n$, the images of $A_k\rightarrow A_r$ and $A_n\rightarrow A_r$ are the same.

\begin{corollary}\label{cor:mf}
    In the situation above, the sequence of reduced twisted $K$-theory groups 
    \begin{equation*}
        \{\widetilde{\KU}^0(EG^n/G)_{\alpha}\}
    \end{equation*}
    satisfies the Mittag-Leffler condition.
    \begin{proof}
        Consider the functor $C:(\text{pro-groups})\rightarrow(\text{topological groups})$ given by taking
        the inverse limit and endowing it with the limit topology and the functor
        $F:(\text{topological groups})\rightarrow(\text{pro-groups})$ given by taking the directed system
        of quotients $G/U$ where $U$ is an open subgroup of $G$. If $\{A_n\}$ is a pro-group, then
        $F(C(\{A_n\}))\iso\{A_n\}$ if and only if $\{A_n\}$ satisfies
        the Mittag-Leffler condition. It is evident that
        \[
            \{\KU_G^0(*)_{\alpha}/I^n\cdot\KU_G^0(*)_{\alpha}\}
        \]
        satisfies the Mittag-Leffler condition because the maps in the system are all
        surjections. By the theorem, it follows that
        \begin{equation*}
            \{\KU_G^0(EG^n)_{\alpha}\}
        \end{equation*}
        does too. But
        \begin{equation*}
            \{\KU_G^0(EG^n)_{\alpha}\}\iso\{\KU^0(EG^n/G)_{\alpha}\}
        \end{equation*}
        since the action of $G$ on $EG^n$ is free. Thus, $\{\KU^0(EG^n/G)_{\alpha}\}$ satisfies the Mittag-Leffler condition.
        Given the fact that the reduced twisted $K$-theory pro-group is the kernel of the rank map
        \begin{equation*}
            \{\KU^0(EG^n/G)_{\alpha}\}\rightarrow\ZZ,
        \end{equation*}
        it follows that the reduced twisted $K$-theory groups satisfies the Mittag-Leffler
        condition.
    \end{proof}
\end{corollary}

We come to our main theorem on twisted equivariant $K$-theory, which we will use in the
next section to construct the examples.

\begin{theorem}\label{thm:groupindices}
    If $G$ is a compact Lie group, then $\ind_{\K}(\alpha)=\ind_{\topo}(\alpha)=\ind_G(\alpha)$ for all
    $\alpha\in\Hoh^3(BG,\ZZ)$. Moreover, $\ind_{\K}(\alpha)$ generates the subgroup
    $\Eoh_{\infty}^{0,0}\subseteq\Eoh_2^{0,0}$.
    \begin{proof}
        It follows from Corollary~\ref{cor:asiso} that there is an isomorphism
        \[
            \KU^0_G(*)_{\alpha}\hat{_I}\iso\KU^0(BG)_{\alpha}.
        \]
        The group $\KU^0_G(*)$ computes the Grothendieck group of $\alpha$-twisted $G$-equivariant vector
        bundles~\cite{freed-hopkins-teleman}*{Proposition~3.5.i} (see also~\cite{freed-hopkins-teleman}*{Example~1.10}),
        and such a bundle is precisely a projective representation of $G$ having obstruction
        class $\alpha$, so there is a natural isomorphism
        \[
            R(G)_{\alpha}\hat{_I}\rightarrow\KU^0(BG)_{\alpha}.
        \]
        All of these groups have compatible rank maps to $\ZZ$.
        As the classes of the augmentation ideal $I$ all have rank $0$,
        the two maps $R(G)_{\alpha}\rightarrow\ZZ$ and $R(G)_{\alpha}\hat{_I}\rightarrow\ZZ$ have
        the same image. Thus, $\ind_{\K}(\alpha)=\ind_G(\alpha)$.

        The second statement follows from Corollary~\ref{cor:mf} and
        Proposition~\ref{prop:kah} applied to the filtration of $BG$ given by $BG^n=EG^n/G$,
        since $\KU^0_G(EG^n)_{\alpha}\iso\KU^0(BG^n)_{\alpha}$.
    \end{proof}
\end{theorem}

\begin{remark}
    When $G$ is a finite group, all the indices in the theorem are finite,
    as one can see by explicitly constructing such representations from cocycles. See~\cite{karpilovsky}.
\end{remark}

\section{Theorems}

% \begin{proposition}
%     Let $X$ be a CW-complex with only finitely many cells in each dimension. If $\ind_{\K}(\alpha)$
%     is finite, then, $\ind_{\K}(\alpha)=\ind(\alpha)$.
%     \begin{proof}
%         Let $X=\colim_k X_k$ be the skeletal filtration of $X$, so that each $X_k$ is a finite
%         CW-complex with cells of dimension at most $k$. Let $\alpha_k=\alpha|_{X_k}$. By a twisted
%         version of the Milnor exact sequence, there is a natural surjection
%         \[
%             \KU^0(X)_{\alpha}\rightarrow\lim_k\KU^0(X_k)_{\alpha_k}\rightarrow 0.
%         \]
%         Suppose that $\ind_{\K}(\alpha)<\ind(\alpha)$.
%         There is some $k$ such that $\ind_{\K}(\alpha_k)=\ind_{\K}(\alpha)$ by the exact
%         sequence. We have $\ind_{\K}(\alpha)=\ind_{\K}(\alpha_m)=\ind(\alpha_m)$ for all $m\geq k$, as $X_m$ is finite.
%         Let $\ind_{\K}(\alpha)=p^ar$ where $p$ is a prime and $r$ is prime to $p$. Then, for all
%         $m\geq k$, there exists a map $X_m\rightarrow\BPU_{p^a s}$ with class $\alpha_m$ for some
%         $s$ prime to $p$ (depending on $m$).
%     \end{proof}
% \end{proposition}

\begin{proposition}
    For any integers $m,n>1$ which have the same prime divisors and such that $m|n$, there
    is a finite abelian group $G$ and a class
    $\alpha\in\Hoh^3(BG,\ZZ)_{\tors}=\Hoh^2(BG,\CC^*)_{\tors}$ such that $\per(\alpha)=m$
    and $\ind_G(\alpha)=n$.
    \begin{proof}
        It suffices by taking products to assume that $m$ and $n$ are powers of the same
        prime. 

        Suppose that $m=p^r$ and $n=p^s=p^{kr+g}$ with $k\geq 1$.
        For any positive integer $t$, let $G_t=(\ZZ/p^t)^2$. Let $G=(G_r)^k\times G_g$. By~\cite{higgs-abelian}*{Lemma~2.2}, for any subgroup $S$ of $G$ such
        that $G/S\iso H\times H$, that is, such that $G/S$ is of
        \textit{symmetric type} in the terminology of \cite{higgs-abelian}, there exists a
        $2$-cocycle $\alpha:G\times G\rightarrow\CC^*$ such that $U(G,\alpha)=S$, where $U(G,\alpha)$ is the group of
        elements $x$ of $G$ such that $\alpha(x,y)=\alpha(y,x)$ for all $y\in G$. Moreover,
        $[G: U(G,\alpha)]=d^2$, where $d$ is the degree of any irreducible projective
        representation with obstruction cocycle $\alpha$.
        Taking $S=\{0\}$ in $G$, we find a $2$-cocycle $\alpha$ such that
        $\ind_G(\alpha)=p^{kr+g}=p^s$.

        It remains to show that $\alpha$
        has the desired period. We know that $\per_{\topo}(\alpha)|p^r$ since
        $G$ is $p^r$-torsion. Let $x$ be an element of $G$ of order $p^r$. Then, $p^{r-1}x$
        is not in $U(G,\alpha)$ since the group is trivial, so there exists $y$
        in $G$ such that $\alpha(p^{r-1}x,y)\neq\alpha(y,p^{r-1}x)$. Let
        $H$ be the subgroup of $G$ generated by $x$ and $y$. Then, $p^{r-1}x$ is not contained
        in $U(H,\alpha_H)$, so that $H/U(H,\alpha_H)$ contains an element
        $\overline{x}$ of order $p^r$. Since the group $H/U(H,\alpha_H)$ is symmetric
        by~\cite{higgs-abelian}*{Lemma~2.2}, it follows that it is isomorphic to $G_r$ and
        that $U(H,\alpha_H)=\{0\}$. Therefore, $\ind_H(\alpha_H)=p^r$, and
        it follows from~\cite{higgs-abelian}*{Lemma~2.3} that this implies that
        $\per_{\topo}(\alpha_H)=p^r$. Thus, $\per_{\topo}(\alpha)=p^r$.
    \end{proof}
\end{proposition}
The next theorem shows that the \'etale index defined
in the first author's thesis is indeed different, in general, from the period.

% As a
% corollary, we find that the indices of unramified Brauer classes in the function fields of
% high dimensional Serre-Godeaux varieties are can be computed purely in terms of projective
% representation theory.

\begin{theorem}\label{thm:main}
    For any integers $m,n>1$ which have the same prime divisors and such that $m|n$, there is a smooth
    projective complex scheme $X$ and a class $\alpha\in\Br(X)$ such that $\per(\alpha)=m$
    and $\eti(\alpha)=\ind(\alpha)=n$.
    \begin{proof}
        Let $G$ be a finite group with a class
        $\overline{\beta}\in\Hoh^3(BG,\ZZ)_{\tors}$ such that $\per_{\topo}(\overline{\beta})=m$
        and $\ind_G(\overline{\beta})=n$. 
        By Theorem~\ref{thm:groupindices}, we have $\ind_{\K}(\overline{\beta})=n$.
        Moreover, by Proposition~\ref{prop:kah}, the group of permanent cycles in
        $\Eoh_2^{0,0}$ is generated by $n$. In particular there are only finitely many
        non-zero differentials leaving $\Eoh^{0,0}$. Let $d_{2s+1}^{\alpha}$ be the last
        non-zero differential leaving $\Eoh^{0,0}$.

%         For example, when $G=(\ZZ/p)^4$,
%         it is shown in~\cite{higgs-healy}*{Theorem~2.2} that of the $p^6-1$ non-zero elements of
%         $\Hoh^3(BG,\ZZ)$, all of which have period $p$, exactly $p^2(p^3-1)(p-1)$ elements have
%         representation index $p^2$.

        We may view the class $\overline \beta$ as a class in $\Br_{\topo}'(BG \times K(\ZZ,2))= \Br_{\topo}'(BG)$,
        and $\per_{\topo}(\overline{\beta})$ and $\ind_{\K}(\overline{\beta})$ are the same computed on either
        $BG$ or $K(\ZZ,2)\times BG$, as one sees by splitting the projection $K(\ZZ,2)\times BG\rightarrow BG$.

        Let $X_{s}$ be a $(2s+2)$-dimensional Serre-Godeaux variety associated to $G$.
        There is a natural morphism $X_s\rightarrow
        K(\ZZ,2)\times BG$, which is a $(2s+1)$-equivalence, and therefore induces an isomorphism in integral cohomology in degrees at most
        $2s$ and an inclusion in degree $2s+1$. Denote the image of $\overline \beta$ in
        $X_{s}$ by $\overline \beta_{s}$; it is non-zero because, for $s \ge 1$, the map
        $\Hoh^3(BG\times K(\ZZ,1), \ZZ) \to \Hoh^3(X_s, \ZZ)$ is an
        inclusion. By direct comparison we find $\ind_\K(\overline \beta_s) | \ind_\K(\overline \beta)$, whereas by our choice of $s$ and
        Proposition~\ref{prop:kah} we find $\ind_\K(\overline \beta) | \ind_\K(\overline \beta_s)$, whereupon it follows that
        $$\ind_{\K}(\overline{\beta_s})=\ind_{\K}(\overline{\beta})=n.$$

        Because the composition
        \[ \Hoh^2_\et(X_s, \mu_m) \iso \Hoh^2(X_s, \ZZ/m) \to \Br'(X_s) \to {}_m\Hoh^3(X_s, \ZZ) \]
        is surjective, we can find a class $\alpha \in \Br'(X)$ whose
        image in $\Hoh^3(X_s, \ZZ)$ is $\overline \beta_n$ and whose period satisfies
        $\per(\alpha) = \per_{\topo}(\overline \beta_s)=m$. The
        relation 
        \[ n=\ind_\K(\overline \beta_s) | \eti(\alpha) | \ind(\alpha)\]
        then follows from~\eqref{eq:1}. This already shows that the \'etale index is
        different, in general, from the period.

        We show that $\eti(\alpha)=\ind(\alpha)=n$ as follows. Write $\mathscr{O}_{X_s}$ for the sheaf of
        holomorphic functions on $X_s$. Since $X_s$ is projective, it is easy to see that
        the natural map $\Br'(X_s)\rightarrow\Hoh^2(X_s,\mathscr{O}_{X_s}^*)_{\tors}$ is
        an isomorphism; see, for example,~\cite{schroer}*{Proposition~1.3}.
        As part of the long exact sequence associated to the exponential sequence
        $0\rightarrow\ZZ\xrightarrow{2\pi i}\mathscr{O}_{X_s}\xrightarrow{\exp}\mathscr{O}_{X_s}^*\rightarrow 1$, we have
        \begin{equation}\label{eq:expex}
            \Hoh^2(X,\mathscr{O}_{X_s})\rightarrow\Hoh^2(X,\mathscr{O}_{X_s}^*)\rightarrow\Hoh^3(X,\ZZ)\rightarrow\Hoh^3(X,\mathscr{O}_{X_s}).
        \end{equation}
        Finally, observe that $\Hoh^2(BG\times K(\ZZ,2),\QQ)\iso\QQ$, which implies that $\Hoh^2(X_s,\mathscr{O}_{X_s})=0$ by the Hodge decomposition.
        Therefore,~\eqref{eq:expex} shows that
        $\Br'(X)\iso\Hoh^2(X,\mathscr{O}_{X_s}^*)_{\tors}\iso\Hoh^3(X,\ZZ)_{\tors}$, since
        $\Hoh^3(X,\mathscr{O}_{X_s})$ is a complex vector space.

        Since $\Br'(X_s)=\Hoh^3(X_s,\ZZ)_{\tors}$, and in particular the comparison map is injective, any flat projective bundle on $X_s$
        with topological Brauer class $\overline{\beta_s}$ has Brauer class $\alpha$. The projective bundles on $X$ which are considered in
        computing $\ind_G(\overline \beta)$ are the bundles arising from a projective representation of $\pi_1(X) = G$, viz.~the flat
        bundles. The restrictions of flat bundles to $X_s$ remain flat, and are consequently algebraizable, and so $\ind(\alpha) |
        \ind_G(\beta_s)=n$, proving that $\eti(\alpha)=\ind(\alpha)=n$, as desired.
    \end{proof}
\end{theorem}

\begin{scholium}
    For any integers $m,n>1$ which have the same prime divisors and such that $m|n$, there is a smooth
    affine complex scheme $X$ and a class $\alpha\in\Br(X)$ such that $\per(\alpha)=m$
    and $\eti(\alpha)=\ind(\alpha)=n$.
    \begin{proof}
        Let $Y$ be a smooth projective scheme $Y$ with a class $\alpha\in\Br(Y)$ with
        $\per(\alpha)=\per_{\topo}(\overline{\alpha})=m$ and
        $\ind_{\topo}(\overline{\alpha})=\eti(\alpha)=\ind(\alpha)=n$, as constructed in the proof of the
        theorem. By Jouanolou's device~\cite{jouanolou}*{Lemme~1.5},
        there exists a vector bundle torsor $\pi:X\rightarrow Y$ with affine total space $X$.
        Topologically, $\pi$ is a homotopy equivalence, so
        \[
        \per(\pi^*\alpha)=\per_{\topo}(\pi^*\overline{\alpha})=m
        \]
        and
        \[
        n=\ind_{\topo}(\pi^*\overline{\alpha})\leq\eti(\pi^*\alpha)\leq \eti(\alpha)=n.
        \]
    \end{proof}
\end{scholium}

In the course of the proof of the theorem, we saw that if $X$ is a Serre-Godeaux variety of
dimension at least $5$, then \[  \Br_{\un}(\CC(X))=\Br(X)=\Hoh^3(X,\ZZ)_{\tors}.  \]
For such a variety and a class $\alpha\in\Hoh^3(X,\ZZ)_{\tors}$, we may thus speak
unambiguously of the algebraic index $\ind(\alpha)$, where we view $\alpha$ as a class in $\Br(X)$.

\begin{theorem}\label{thm:indsg}
    Let $G$ be a finite group. There exists a positive integer $n$ such that if $X$ is a
    Serre-Godeaux variety for $G$ of dimension at least $n$, and if
    $\alpha\in\Hoh^3(BG,\ZZ)$ is a Brauer class, then $\ind(\alpha_{\CC(X)})=\ind(\alpha)=\ind_G(\alpha)$. That is, the index is
    computed by projective representations.
    \begin{proof}
        The theorem follows from the proof of Theorem~\ref{thm:main}.
        We can find a single finite $n$ because $\Hoh^3(BG,\ZZ)$ is a finite abelian group.
        This shows that we can compute $\ind(\alpha)=\ind_G(\alpha)$ for $\alpha\in\Br(X)$. It follows from
        an argument of David Saltman that $\ind(\alpha_{\CC(X)})=\ind(\alpha)$, since $X$ is
        regular and noetherian; see~\cite{aw3}*{Proposition~5.5}.
    \end{proof}
\end{theorem}

In the theorem, it is clear that $\ind(\alpha)|\ind_G(\alpha)$. The role of the \'etale
index or the topological index
is to show that the index is large enough, and that this divisibility relation is actually
an equality.


\begin{bibdiv}
\begin{biblist}

\bib{antieau-cech}{article}{
    author={Antieau, Benjamin},
    title={\v Cech approximation to the Brown-Gersten spectral sequence},
    journal={Homology Homotopy Appl.},
    volume={13},
    date={2011},
    number={1},
    pages={319--348},
    issn={1532-0073},
}

\bib{antieau}{article}{
    author = {Antieau, Benjamin},
    title = {Cohomological obstruction theory for Brauer classes and the period-index problem},
    journal = {Journal of K-theory},
    volume = {8},
    number = {3},
    pages={419--435},
    date = {2011},
}

\bib{aw1}{article}{
    author = {Antieau, Benjamin},
    author = {Williams, Ben},
    title = {The period-index problem for twisted topological K-theory},
    journal = {ArXiv e-prints},
    eprint = {1104.4654},
    year = {2011},
}

\bib{aw3}{article}{
    author = {Antieau, Benjamin},
    author = {Williams, Ben},
    title = {Low-dimensional topological period-index problems},
    year = {2012},
}

\bib{atiyah-hirzebruch}{article}{
    author={Atiyah, Michael F.},
    author={Hirzebruch, F.},
    title={Analytic cycles on complex manifolds},
    journal={Topology},
    volume={1},
    date={1962},
    pages={25--45},
    issn={0040-9383},
}

\bib{atiyah-segal}{article}{
    author={Atiyah, Michael F.},
    author={Segal, Graeme},
    title={Twisted $K$-theory},
    journal={Ukr. Mat. Visn.},
    volume={1},
    date={2004},
    number={3},
    pages={287--330},
    issn={1810-3200},
    translation={
        journal={Ukr. Math. Bull.},
        volume={1},
        date={2004},
        number={3},
        pages={291--334},
        issn={1812-3309},
    },
}

\bib{atiyah-segal-cohomology}{article}{
    author={Atiyah, Michael F.},
    author={Segal, Graeme},
    title={Twisted $K$-theory and cohomology},
    conference={
        title={Inspired by S. S. Chern},
    },
    book={
        series={Nankai Tracts Math.},
        volume={11},
        publisher={World Sci. Publ., Hackensack, NJ},
    },
    date={2006},
    pages={5--43},
}

% \bib{baumbrowder}{article}{
% 	title = {The cohomology of quotients of classical groups},
% 	volume = {3},
% 	journal = {Topology},
% 	author = {Baum, Paul F.},
%     author={Browder, William},
% 	year = {1965},
% 	pages = {305--336}
% }

\bib{boardman}{incollection}{
    AUTHOR = {Boardman, J. Michael},
     TITLE = {Conditionally convergent spectral sequences},
 BOOKTITLE = {Homotopy invariant algebraic structures ({B}altimore, {MD}, 1998)},
    SERIES = {Contemp. Math.},
    VOLUME = {239},
     PAGES = {49--84},
 PUBLISHER = {Amer. Math. Soc.},
   ADDRESS = {Providence, RI},
      YEAR = {1999},
}

% \bib{dejong}{article}{
% 	title = {{A} result of {Gabber}},
% 	author = {A. J. de Jong},
% 	year = {2003},
% 	note = {unpublished preprint, \texttt{www.math.columbia.edu/{\textasciitilde}dejong/papers/2-gabber.pdf}}
% }

% \bib{dwyer}{article}{
%     author={Dwyer, Christopher},
%     title={Twisted equivariant $K$-theory for proper actions of discrete
%     groups},
%     journal={$K$-Theory},
%     volume={38},
%     date={2008},
%     number={2},
%     pages={95--111},
%     issn={0920-3036},
% }

\bib{freed-hopkins-teleman}{article}{
    author={Freed, Daniel S.},
    author={Hopkins, Michael J.},
    author={Teleman, Constantin},
    title={Loop groups and twisted $K$-theory I},
    journal={J. Topol.},
    volume={4},
    date={2011},
    number={4},
    pages={737--798},
    issn={1753-8416},
}

\bib{grothendieck-brauer}{article}{
    author={Grothendieck, Alexander},
    title={Le groupe de Brauer. I. Alg\`ebres d'Azumaya et interpr\'etations diverses},
%     language={French},
    conference={
    title={S\'eminaire Bourbaki, Vol.\ 9},
    },
    book={
    publisher={Soc. Math. France},
    place={Paris},
    },
    date={1995},
    pages={Exp.\ No.\ 290, 199--219},
}

\bib{higgs}{article}{
    author={Higgs, R. J.},
    title={On the degrees of projective representations},
    journal={Glasgow Math. J.},
    volume={30},
    date={1988},
    number={2},
    pages={133--135},
    issn={0017-0895},
}

\bib{higgs-abelian}{article}{
    author={Higgs, R. J.},
    title={Projective representations of abelian groups},
    journal={J. Algebra},
    volume={242},
    date={2001},
    number={2},
    pages={769--781},
    issn={0021-8693},
%     review={\MR{1848971 (2002f:20014)}},
%     doi={10.1006/jabr.2000.8751},
}

% \bib{higgs-healy}{article}{
%     author={Higgs, Russell},
%     author={Healy, Donal},
%     title={Projective character degree patterns of groups of order $p^4$},
%     journal={Comm. Algebra},
%     volume={34},
%     date={2006},
%     number={12},
%     pages={4623--4630},
%     issn={0092-7872},
% %     review={\MR{2273729 (2007k:20024)}},
% %     doi={10.1080/00927870600936831},
% }

\bib{jouanolou}{article}{
    author={Jouanolou, J. P.},
    title={Une suite exacte de Mayer-Vietoris en $K$-th\'eorie alg\'ebrique},
%     language={French},
    conference={
        title={Algebraic $K$-theory, I: Higher $K$-theories (Proc. Conf., Battelle Memorial Inst., Seattle, Wash., 1972)}, },
        book={
            publisher={Springer},
            place={Berlin},
    },
    date={1973},
    pages={293--316. Lecture Notes in Math., Vol. 341},
}

\bib{karpilovsky}{book}{
    author={Karpilovsky, Gregory},
    title={Projective representations of finite groups},
    series={Monographs and Textbooks in Pure and Applied Mathematics},
    volume={94},
    publisher={Marcel Dekker Inc.},
    place={New York},
    date={1985},
    pages={xiii+644},
    isbn={0-8247-7313-6},
}

\bib{lahtinen}{article}{
    author = {Lahtinen, Anssi},
    title = {The Atiyah--Segal completion theorem in twisted K-theory},
    journal = {ArXiv e-prints},
    eprint = {0809.1273},
    year = {2008},
}

\bib{milnor-construction}{article}{
    author={Milnor, John},
    title={Construction of universal bundles. II},
    journal={Ann. of Math. (2)},
    volume={63},
    date={1956},
    pages={430--436},
    issn={0003-486X},
}

\bib{milnor}{article}{
    author={Milnor, John},
    title={On axiomatic homology theory},
    journal={Pacific J. Math.},
    volume={12},
    date={1962},
    pages={337--341},
    issn={0030-8730},
}


\bib{schroer}{article}{
    AUTHOR = {Schr{\"o}er, Stefan},
    TITLE = {Topological methods for complex-analytic {B}rauer groups},
    JOURNAL = {Topology},
    VOLUME = {44},
    YEAR = {2005},
    NUMBER = {5},
    PAGES = {875--894},
    ISSN = {0040-9383},
}

\bib{serre}{article}{
    author={Serre, Jean-Pierre},
    title={Sur la topologie des vari\'et\'es alg\'ebriques en caract\'eristique $p$},
    conference={
        title={Symposium internacional de topolog\'\i a algebraica},
    },
    book={
        publisher={Universidad Nacional Aut\'onoma de M\'exico and UNESCO, Mexico City},
    },
    date={1958},
    pages={24--53},
}




\end{biblist}
\end{bibdiv}
\end{document}